\documentclass{amsart}
\usepackage{changebar}
\usepackage{lineno}
\usepackage{amsthm}
\usepackage{verbatim}
\usepackage[utf8]{inputenc}
\usepackage{blindtext}

\makeatletter
\renewcommand\part{%
   \if@noskipsec \leavevmode \fi
   \par
   \addvspace{4ex}%
   \@afterindentfalse
   \secdef\@part\@spart}

\def\@part[#1]#2{%
    \ifnum \c@secnumdepth >\m@ne
      \refstepcounter{part}%
      \addcontentsline{toc}{part}{\thepart\hspace{1em}#1}%
    \else
      \addcontentsline{toc}{part}{#1}%
    \fi
    {\parindent \z@ \raggedright
     \interlinepenalty \@M
     \normalfont
     \ifnum \c@secnumdepth >\m@ne
       \Large\bfseries \partname\nobreakspace\thepart
       \par\nobreak
     \fi
     \huge \bfseries #2%
     \par}%
    \nobreak
    \vskip 3ex
    \@afterheading}
\def\@spart#1{%
    {\parindent \z@ \raggedright
     \interlinepenalty \@M
     \normalfont
     \huge \bfseries #1\par}%
     \nobreak
     \vskip 3ex
     \@afterheading}
\makeatother

\newtheorem{lemma}{Lemma}
\newtheorem{theorem}{Theorem}
\newtheorem{proposition}{Proposition}
\newtheorem{corollary}{Corollary}
\newtheorem{definition}{Definition}

\newtheorem{assumption}{Assumption}
\newtheorem{remark}{Remark}
\newtheorem{fact}{Fact}

\newcommand{\nth}[1]{$#1 {\rm - th }$}

\newcommand{\Z}{\mathbb{Z}}
\newcommand{\ZM}[1]{\Z /( #1 \cdot \Z)}
\newcommand{\ZMs}[1]{(\Z / #1 \cdot \Z)^*}

\newcommand{\ord}{{\rm ord}}

\newcommand{\Tr}{\mbox{\bf Tr}}

\newcommand{\rg}[1]{\mbox{\bf #1}}
\newcommand{\eu}[1]{\mathfrak{#1}}
\newcommand{\id}[1]{\mathcal{#1}}
\newcommand{\Gal}{\mbox{ Gal }}

\newcommand{\Ker}{\mbox{ Ker }}

\newcommand{\rf}[1]{(\ref{#1})}

\newcommand{\nmid}{\not \hspace{0.25em} \mid}
\newcommand{\Norm}{\mbox{\bf N}}

\newcommand{\F}{\mathbb{F}}

\newcommand{\K}{\mathbb{K}}

\newcommand{\N}{\mathbb{N}}
\newcommand{\Q}{\mathbb{Q}}

\def\ra{\rightarrow}
\newcommand{\ran}{\rangle}
\newcommand{\lan}{\langle}

\newcommand{\ol}{\overline}
\newcommand{\veps}{\varepsilon}

\newcommand{\st}{^{\times}}

\begin{document}
 {\obeylines \small
 \vspace*{-0.2cm}
 \hspace*{3.5cm}Par la m\`ere apprenant que son fils est gu\'eri,
\hspace*{3.5cm}par l'oiseau rappelant l'oiseau tomb\'e du nid,
\hspace*{3.5cm}par l'herbe qui a soif et recueille l'ond\'ee,
\hspace*{3.5cm}par le baiser perdu par l'amour redonn\'e,
\hspace*{3.5cm}et par le mendiant retrouvant sa monnaie:
\vspace*{0.2cm} 
\hspace*{4.0cm} Je vous salue, Marie \footnote{Francis James: {\em Pri\`ere}. Music by Georges Brassens}
\vspace*{0.4cm}
\hspace*{4.5cm}\indent    {\it To the memory of Maria, Marta and Montserrat;}
\hspace*{5.5cm}\indent    {\it to our wives and children}
\vspace*{0.4cm}
\smallskip}

%{\obeylines \small
%\vspace*{-0.2cm}
%\hspace*{3.5cm}Par la m\`ere apprenant que son fils est gu\'eri,
%\hspace*{3.5cm}par l'oiseau rappelant l'oiseau tomb\'e du nid,
%\hspace*{3.5cm}par l'herbe qui a soif et recueille l'ond\'ee,
%\hspace*{3.5cm}par le baiser perdu par l'amour redonn\'e,
%\hspace*{3.5cm}et par le mendiant retrouvant sa monnaie:
%\vspace*{0.2cm} 
%\hspace*{4.0cm} Je vous salue, Marie \footnote{Francis James: {\em Pri\`ere}. Music by Georges Brassens}
%\vspace*{0.4cm}
%\hspace*{4.5cm}\indent    {\it To the memory of Maria, Marta and Montserrat;}
%\hspace*{5.5cm}\indent    {\it to our wives and children}
%\vspace*{0.5cm}
%\smallskip}
\title[Beyond Fermat desert ] {Beyond and slightly aside of the Fermat desert} 
\author{Boris Bartolom\'{e} and Preda Mih\u{a}ilescu} 
\address[P. Mih\u{a}ilescu]{Mathematisches Institut der Universit\"at 
G\"ottingen}
\email[P. Mih\u{a}ilescu]{Preda@uni-math.gwdg.de} 
\address[B. Bartolom\'{e}]{}
\email[B. Bartolom\'{e}]{Boris.Bartolome@mathematik.uni-goettingen.de} 
\date{Version 1.0 \today}
\vspace{0.3cm}
%\huge
\begin{abstract}
We give a cyclotomic proof of the fact that the equation $\frac{x^p + y^p}{x+y}  =  p^e z^p$ has no 
solutions in coprime integers $x,y,z$ and $p > 3$, a prime. This implies in particular Fermat's Last Theorem.
\end{abstract} 
\maketitle
\vspace*{-0.5cm}
\tableofcontents
%\vspace*{0.5cm}
%\newpage
\vspace*{1.0cm}
\section{Introduction}
The reader familiar with the amazing diversity and richness of criteria developed until $1990$, in the attempt
to prove that the Fermat Conjecture, purporting that the equation
\begin{eqnarray}
    \label{flt}
    x^p + y^p + z^p = 0; \quad x, y, z \in \Z; \quad x y z \neq 0, \quad p > 2
\end{eqnarray}

has no solutions, can observe that it long resembled a vast and barren desert, sporadically punctuated by rare and isolated oasis—exceptional, innovative, and surprising conditions that would necessarily have to hold if the conjecture were false. These conditions were so intricate and specific that their fulfillment appeared highly implausible. Nevertheless, none of these investigative approaches succeeded in reaching the definitive resolution—the metaphorical “ocean” beyond the desert—where the ultimate truth of the conjecture would be established.
This description is particularly apt in the context of the {\em First Case}, which assumes that the prime $p \nmid x y z$. Ribenboim’s extensive and illuminating survey \cite{Ri} captures the breadth of attempts made, documenting the numerous avenues pursued and the audacity of those efforts. 
\subsection{Recent history}
(To be reinserted after complete referee's review of the paper)

\subsection{Overview of the paper}
In this paper we consider the following Diophantine norm equation: 
\begin{eqnarray}
\label{sft}
\quad \frac{x^p + y^p}{x+y} & = & p^e z^p; \quad (x, y, z) \in (\Z^*)^3, \quad ( x,y,z ) = 1, \\ \nonumber & &  
\hbox{and $p > 3$ is a prime; $e \in \{0, 1\}$ }.
\end{eqnarray}
We will refer to the case when $e=0$ as the first case and, when $e = 1$, as the second case of equation \rf{sft}.
The relation between the value of $e$ and $x, y, z$ is explained below. 
The connection to \rf{flt} will become apparent below, from the classical formulas of Barlow and Abel. 
The term {\em Strong Fermat, or SFLT} has been coined by Georges Gras and Roland Qu\^eme \cite{GQ} for equation \rf{sft}. 
We shall use, by analogy, the term {\em Strong Fermat-Catalan, or SFC} for the modified version of the equation 
\rf{sft}, in which the exponent of $z$ is an arbitrary prime $q$:
\begin{eqnarray}
\label{sfct}
\quad \frac{x^p + y^p}{x+y} & = & p^e z^q; \quad (x, y, z) \in (\Z^*)^3, \quad ( x,y,z ) = 1, \\ \nonumber & &  
\hbox{and $p > 3$ is a prime; $e \in \{0, 1\}$ }.
\end{eqnarray}
We prove:
\begin{theorem}
    \label{main}
    The Diophantine equation \rf{sft} has no integer solutions.
\end{theorem}
The method of proof is elastic enough, as to apply with slight modifications also for \rf{sfct}, a fact which 
we prove in a paper prepared in parallel to the present one. It follows from Barlow and Abel (\cite{Ri}, Lecture IV), that the
Fermat Equation can be considered as an equation \rf{sft} with additional constraints; the truth of the
Fermat Conjecture follows thus a fortiori from Theorem \ref{main}. We note that for $p=3$, \rf{sft} 
has infinitely many solutions that stem from norms of \nth{p} powers in the Eisenstein integers: 
if $(x, y) \in \Z^2$, $\rho$ is a complex third root of unity and 
$\xi^3 = x + y \rho$ for some $\xi \in \Z[\rho]$, then $\frac{x^3+y^3}{x+y} = z^3$, with $z = \xi \cdot \bar{\xi} \in \Z^*$.
The condition $p > 3$ is thus necessary; recall also that Fermat's Conjecture had been already proved by Euler
for $p = 3$.

The two proofs of Fermat's Conjecture discussed above, were made possible by the proof of some independent
and very powerful general theorems, with wide consequences in different areas of mathematics, and which just 
{\em so happened} to be strong enough to {\em also} imply Fermat's Conjecture. Our approach is situated at the other
end of the spectrum: we use methods restricted to the class field theory of the $p$-cyclotomic field, methods
which essentially were known already more than $100$ years ago, and our proof appears as a remote echo to some
of the results derived at the beginning of the $20$-th century. Explicitly, Furtw\"angler proved the following
statement: if $x, y, z$ is a solution to \rf{flt} with $p \nmid x y z$, then the primes $r | x y z$ verify
$r^{p-1} \equiv 1 \bmod p^2$. Equivalently, the prime ideals $\eu{R} \subset \Q[ \zeta_p ]$ laying above
$r$ are further totally split in the \nth{p^2} cyclotomic extension \cite{Fur}. With this, Furtw\"angler could derive
results previously discovered by Mirimanoff and Wieferich, and offer some generalizations thereof. For his
proof, Furtw\"angler applied class field theory exclusively by use of Eisenstein's $p$-power explicit reciprocity
laws. Interestingly, using class group annihilation by the Stickelberger ideal, Furtw\"angler's result receives
a shorter and more transparent proof. The interested reader will be provided, in the next chapter, basic forumulae
used to reduce the proof of Furtw\"angler reduces to a few lines. 

In our proof, we make intensive use of a large submodule $J$ of the Stickelberger ideal $I$, which 
comes with pleasant additional properties with respect to possible solutions of \rf{sft}. Explicitly, for
$t \in J$ and $(x, y, z)$ a solution of \rf{sft}, there is a $\beta( t ) \in \Z[ \zeta_p ]$ with
\[ \beta( t )^p = \frac{(y + \zeta_p x )^t}{p^{ek}}, \quad k \geq 1.\]
The numerical discriminant of $\beta( t )$ induces an intricate factorization of the invariant $K = x y (x + y)$, which is coprime to $z$, by Remark \ref{remark1} below. A detailed investigation of this factorization leads to the proof of our fundamental Proposition \ref{pmain}, 
which states that the rational primes $r | K$ are totally split in the real $p$-cyclotomic 
field $\Q[ \zeta_p + \bar{\zeta}_p ]$, thus verifying $r \equiv \pm 1 \bmod p$. One easily checks that
$K$ must be even, and so $2$ should verify the same congruence, which is impossible. The formal resemblance with
Furtw\"angler's result is quite impressive -- given that the methods of proof do not share much in common. 

Given that $z$ does not occur in our criterion, the proof holds beyond the particular Fermat Equation. The
classical equation of Nagell-Ljunggren \cite{RiC} can presently also be reduced to the statement that
\begin{eqnarray}
\label{nl}
    \frac{x^p-1}{x-1} = p^e y^q, \quad (x, y) \in (\Z^*)^2, \quad e \in \{0, 1\},
\end{eqnarray}
has no solutions for odd primes $p, q$. The equations relates by factorizations similar to Barlow-Abel, 
to the classical equation of Catalan $x^p - y^q = 1$, the relation being analogous to the one between
\rf{flt} and Strong Fermat. While the study of \rf{nl} began in the early $20^{th}$ century, long before
the first important progress made on Catalan by use of Baker's linear forms in logarithms, the study of
\rf{sft} as an equation of interest in itself was only proposed by Gras and Qu\^eme some time after
Wiles's proof of Fermat. To this day, little is known about this equation and the main known results 
are due to the fact that the most powerful classical criteria for the First Case of Fermat's Last Theorem, 
the criteria of Eichler
\begin{comment}
    \footnote{If $p$ is an odd prime, $x, y, z$ are relative prime integers such that $x^p+y^p+z^p=0$ and $r \in \N$ such that $r|(x-y), p \nmid (x^2-y^2)$, then $r \equiv 1 [p^2]$}
\end{comment}
 \cite{Ei} and (Vandiver-) Sitaraman
 \begin{comment}
 \footnote{For every $n \geq 1$, the equation $x^l^n+y^l^n+z^l^n=0$ has no integral solution $(x,y,z)$ with $l \nmid xyz$ when combined with \cite{Ku}}    
 \end{comment}
 \cite{Sit}, respectively, can be derived using only properties
of \rf{sft}, by imposing the condition $p \nmid K$. Under this condition -- which we refer to
as {\em First Case of \rf{sft}} -- the two conditions previously quoted (Eichler's and Sitaraman's) need to be fulfilled by possible solutions of 
\rf{sft}.

It is interesting to observe that none of the previous proofs of FLT can extend to the case of \rf{sft}.
While a Frey curve for \rf{sft} can be defined, the Ribet descent fails to provide a contradiction and
produces an obstruction to the generalization of Wiles's result to \rf{sft}. In the case of $abc$-inequalities,
the application to \rf{sft} is made challenging by the lack of lower bounds to possible solutions; in fact
the only lower bound known to the authors, stems from the fact that in \rf{sft}, the number $z$ is totally
split in the \nth{p} cyclotomic field, and thus $| z | \geq 2p+1$. The situation is not unsimilar to the
applications of Baker's theory for the proof of Catalan's equation. In all these cases, some extremely powerful
and general results in algebraic geometry -- for Fermat -- and Diophantine approximation, for Catalan, succeed
to reveal valuable information, resp. provide even proofs of the respective classical Conjectures on Diophantine 
Equations. It is however their very generality which sets limitations to the extension of their use to
larger classes of Diophantine Equations. In our case, {\em after finding an adequate entry point} -- which 
is apparently the main challenge -- using methods from the field of cyclotomy, in which the questions are defined, 
one obtains proofs with some more degrees of freedom, that allow for slight generalizations. But these reduce
to the strict domain of {\em Diophantine cyclotomic norm equations} and do not apply outside of
this relatively restricted area of mathematics. 

Andrew Beal investigated solutions of the more general equation
\begin{eqnarray}
\label{beal}
x^p + y^q + z^r = 0
\end{eqnarray}
with prime exponents $p, q, r$ and integers $x, y, z$ and found, by astutely leading computer
searches, that the equation does have solutions when $(x, y, z) > 1$. He then formulated the 
conjecture that \rf{beal} has no solutions in mutually coprime $(x, y, z)$ and exponents
$p, q, r \geq 3$. This conjecture bears the name of {\em Beal Conjecture}. The slightly more general
equation, in which the exponents are allowed to be $2$ or $4$ is called {\em Super Fermat Equation},
and a list of sporadic solutions with some even exponents is known and conjectured to be complete.
The Fermat - Catalan equation $x^p + y^p = z^q$ -- with odd primes $p \neq q$ -- appears thus as a particular 
case of \rf{beal}. The mention we make here of \rf{beal} intends to prevent natural questions that some of
our readers may ask, as to possible extensions of our approach to \rf{beal}. It is said to be unwise to 
state that some result is {\em impossible}. We may, however, safely state the major obstruction which
calls to prudence -- or even strong skepticism -- with respect to a possible extension of our methods to \rf{beal}.
While \rf{sft}, and also \rf{sfct} are strongly connected to ideals of the \nth{p} cyclotomic field -- as shown in 
the following chapter --, allowing the use of the very well understood class field theory based on Stickelberger 
ideals, we need \rf{beal} to work in larger metabelian extensions of cyclotomic fields. Beyond basic representation
theory, no useful explicit results on class field theory are known to us, that might offer themselves to build
a bridge between our present kind of proof and a wishful approach to \rf{beal}.

On the other hand, together with the forthcoming result on \rf{sfct}, our approach may be more successful for 
the investigation of similar equations, defined over some number fields. The authors will probably discontinue
their own work on Diophantine Equations here, but we are looking forward to possible research by colleagues
and their students, aiming to choose well balanced families of number fields in which 
to test the application of our approach for investigating the generalizations of Fermat and Fermat-Catalan
equations. 

Beyond this introduction, the paper is structured in three chapters. In the first, we recall basic
and classical facts about \rf{flt} and \rf{sft}, the Stickelberger Ideal\footnote{While we use in this paper the term of 
Stickelberger ideal for the module of annihilators of the minus part of the class group in cyclotomic 
fields -- see \cite{Wa}, \S 15.1 -- it is interesting to recall that the explicit annihilators for the case of
the $p$-cyclotomic fields, which are the only ones needed here, were already known to Kummer. He did however use less
 their structure as a submodule of the group ring, which is of good use in our proofs.} and its application to
possible solutions of \rf{sft}. \begin{comment}
    The third chapter is an intermezzo on semilocal rings of completions
at primes dividing $K$ over $\Q[ \zeta_p ]$. The rationale for this will be discussed in site, suffice
to mention here the fact that these larger structure offer a perfect framework to organize the ideas of the
proof. The proof itself can however be derived by simple congruences modulo ideals of $\Q[ \zeta_p ]$. In the
last chapter, we introduce definitions aiming at facilitating the investigation of the common factors of
$K$ and the discriminant of $\beta( t )$ described above.
\end{comment}
This leads to the proof of  Proposition \ref{pmain} and of the main Theorem. 

\section{Classical results, notations and prerequisites}
A classical fact, often attributed to Euler, states that 
for coprime integers $x, y$ and $n \in 2\N+1,\quad n>1$, the greatest common divisor
$d = \left( \frac{x^n+y^n}{x+y}, x+y \right)$ divides $n$. This can 
be easily verified by using the substitution $s = x+y$ and by
introducing it in the fraction $\frac{x^n+y^n}{x+y} = \frac{x^n + 
(s-x)^n}{s}$: the  claim follows from $( x, y ) = 1$. In our case, $n = p$ and if $d = 1$, then 
$p \nmid z$, while for $d = p$, the same substitution implies
that $v_p\left(\frac{x^p + y^p}{x+y}\right) = 1$, hence the introduction of $p^e$ in \rf{sft}.
% -- see e.g. \cite{Mi2}. 

Using the methods of Barlow and Abel, we find the generalization of their classical factorization (\cite{Ri}, lecture IV):

There exist integers $u$ and $v$, $(u,v) = 1$ and $p \nmid uv$, such that:

\begin{eqnarray}
\label{factf}
x + y  =  u^p &\hbox{and} & \quad \frac{x^p + y^p}{x+y} = v^p \quad \hbox{if $p \nmid (x+y)$, and} \\
\nonumber x+y  =  p^{p-1} \cdot u^p &\hbox{and} & \quad \frac{x^p + y^p}{x+y} = p v^p \quad \hbox{otherwise}.
\end{eqnarray}
We see that having solutions to the Strong Fermat equation \rf{sft} is a necessary -- but not sufficient -- 
condition for the Fermat equation \rf{flt} to have solutions. Thus, results on solutions of \rf{sft}
imply the same conditions for solutions of \rf{flt}, while the homogenous nature of the latter will 
provide additional information, in form of strong lower bounds on the solutions -- see \cite{Ri}, \textbf{XII}.

\begin{remark}
\label{remark1}
Since $(x+y,\frac{x^p+y^p}{x+y})=1$, it follows that $1=(z^p,x+y) = (z,x+y)$, so $z$ and $xy(x+y)$ are coprime.
\end{remark}

\subsection{Notations}
Throughout this paper, for $r$ a prime or a prime power, we denote
by $\F_r$ the field with $r$ elements. 

We let $P = \{0, 1, \ldots, p-1\}, P^* = P \setminus \{ 0 \}$ be the minimal positive representatives 
for $\F_p$ and $\F_p^{\times}$, respectively. $\zeta$ will be a primitive \nth{p} root of unity and $\K = \Q[ \zeta ]$ the \nth{p} 
cyclotomic field, with Galois group $G = \Gal( \K /\Q )$. The automorphisms $\sigma_c \in G$ are given by
 $\zeta \mapsto \zeta^c$, for $c \in P^*$. We fix $\sigma \in G$, an automorphism generating $G$ as a cyclic group. The Teichm\"uller character is 
 \[ \varpi : G \ra \Z_p\st ; \quad \sigma_c \mapsto \varprojlim_n (c^{p^n} \bmod p^{n+1}) ,\quad c \in P^*. \]
The character restricts to characters mapping $G$ to $\F_p$ and also $\ZM{p^n}$,
and it can also be interpreted as a character of $\F_p\st$ 
via the isomorphism $G \cong \F_p\st$.
It induces the spectral decomposition of $\rg{R}[ G ]$, for $\rg{R} \in \{ \F_p, \Z_p \}$
 as follows: 
 \begin{eqnarray}
 \label{idemps}
 \veps_k & := & \frac{1}{p-1} \sum_{c \in P^*} \varpi^k( c ) \sigma_c^{-1}; \ k \in P^*, \quad \hbox{verifying:} \nonumber \\
 \sum_{k \in P^*} \veps_k & = & 1; \quad \quad \veps_k \cdot ( \sigma_c - \varpi( c )^k ) = 0.
 \end{eqnarray}
 
 \begin{assumption}
 \label{sol}
 We assume in the sequel that $(x, y, z) \in \Z^3$ is a solution to \rf{sft}, for the prime exponent $p$, and $ x  > | y | \geq 1$. 
 We endow $\K$ with the normal base $\id{E} = \{ \zeta^c \ : \ c \in P^*\}$, as a $\Q$-vector space.
\end{assumption}

 The complex conjugation acting on $\K$ is denoted by 
 $\jmath = \sigma_{p-1} = \sigma^{(p-1)/2}$. 
 We use the uniformizor $\lambda = 1-\zeta \in \Z[ \zeta ]$,  that generates the principal prime 
 $\wp \subset \Z[ \zeta ]$ above $p$ (we remind that $\wp^{p-1} = (\lambda)^{p-1} = (p)$).  It induces $\lambda$-adic expansions of algebraic integers in $\Z[ \zeta ]$,
 such that, for $w \in \Z[\zeta ]$, we may write: 
 \[ w = \sum_{j=0}^{\infty} a_j \lambda^j = a_0 + a_1 \lambda + a_2 \lambda^2 + O( \lambda^3 ), \]
where the $a_j \in \{ -\frac{p-1}{2}, -\frac{p-3}{2}, \ldots, \frac{p-1}{2} \}$, 
and only finitely many coefficients are non-null. The symbol $O( \lambda^k )$ stands 
for a remainder term divisible by $\lambda^k$. The same notation can be used also in $\Z_p[ \zeta ]$; thus, for $\alpha \in  \Z_p[ \zeta ]$,
\[ \alpha = \sum_{j=0}^{\infty} a_j \lambda^j . \]

\subsection{Characteristic numbers and characteristic ideals}
In the cyclotomic field $\K$, the norm decomposes as:
\[ \frac{x^p + y^p}{p^{e}(x+y)}  = \prod_{c \in P^*} \frac{y + \zeta^c x}{(1-\zeta^c)^e}. \]
This leads naturally to the following 
\begin{definition}
\label{charids}
We define
\begin{eqnarray}
\label{charnu}
 \alpha = \frac{y + \zeta x}{(1-\zeta)^e},
\end{eqnarray}
as the {\em \textbf{characteristic number}} of equation \rf{sft} and
\begin{eqnarray}
\label{charid}
\eu{A} = ( \alpha, z ) \subset \Z[ \zeta ],
\end{eqnarray}
as the {\em \textbf{the characteristic ideal}} of the same equation.
\end{definition}
One notices that $\Norm_{\K/\Q}( \alpha ) = z^p$, so the characteristic number
encodes the existence of a non-trivial solution to \rf{sft}. By consequence of
the definition, so does the characteristic ideal $\eu{A}$. The following Lemma
gathers important properties of the characteristic number and ideal. 
 
\begin{lemma}
\label{aux1}
\begin{itemize}
\item[ 1. ] The characteristic number $\alpha$ is integral. 
\item[ 2. ] The Galois group $G$ acts on the characteristic number, giving rise to pairwise coprime integral elements; that is,  for $1 \leq c < d \leq p-1$,
\begin{eqnarray*}
( \sigma_c( \alpha ), \sigma_d( \alpha )) = (1) .
\end{eqnarray*}
\item[ 3. ] $\eu{A}$ is related to the {\em characteristic number} $\alpha$ by the relations:
\begin{eqnarray}
\label{rclass}
\eu{A}^p & = & ( \alpha ), \quad \Norm( \eu{A} ) = (z).
\end{eqnarray}
\item[ 4. ] The characteristic number satisfies:
\begin{eqnarray}
\label{al'}
\frac{\alpha}{\bar{\alpha}} & = & \upsilon \cdot \frac{1 + \zeta ( x/y )}{1 + \bar{\zeta}(x/y)}, \quad \hbox{with } \\ \nonumber
\upsilon & = & \begin{cases}
      1 & \hbox{if $ e = 0$, and } \\
      -\bar{\zeta} & \hbox{for $e = 1$,}
    \end{cases}
\end{eqnarray}
and if $e = 0$ we have 
\begin{eqnarray}
\label{twist}
\alpha' = \zeta^{x/(x+y)} \alpha = (x+y)\cdot ( 1 + O(\lambda^2 )).
\end{eqnarray}
\end{itemize} 
\end{lemma}
\begin{proof}
Let's first prove point 1: for $e=0$,  $\alpha$ is clearly an integral element, while if $e=1$, then $p | ( x+y)$, so $\alpha$ is also integral.

For point 2., let's first remind that $\lambda = 1 - \zeta$ and that, for distinct $c, d \in P$, we have $\frac{\zeta^c - \zeta^d}{\lambda}$ is a unit in $\Z[\zeta]$. Let $I( c,d ) = \left( \sigma_c( \alpha), \sigma_d( \alpha )\right) $; then $y \lambda \in I( c, d)$. If $e= 0$,
this follows from $\sigma_c( \alpha ) - \sigma_d( \alpha ) = ( \zeta^c -\zeta^d ) y \in I( c, d)$, and for $e= 1$,
we have $(1-\zeta^c) \sigma_c(\alpha) - (1-\zeta^d) \sigma_d(\alpha) = -( \zeta^c -\zeta^d ) y \in I( c, d)$. 
Likewise, $x \lambda \in I( c, d)$: for $e= 0$,
we have $\bar{\zeta}^c \sigma_c( \alpha ) - \bar{\zeta}^d \sigma_d( \alpha ) = ( \bar{\zeta}^c -\bar{\zeta}^d ) x \in I( c, d)$, 
while for $e= 1$,
we have $(1-\bar{\zeta}^c) \sigma_c( \alpha ) - (1-\bar{\zeta}^d) \sigma_d( \alpha ) =( \bar{\zeta}^d -\bar{\zeta}^c) x \in I( c, d)$. 
Recall that $\lambda = 1-\zeta$ and $\frac{\zeta^a - \zeta^b}{\lambda} \in \id{O}^{\times}( \K )$ for any distinct  $a, b \in P$. Thus
$I( c,d ) | ( x,y) ( \lambda ) = \wp$, since $(x,y) = 1$. However, $( \alpha, p) = (1)$ by definition, so it follows that $I( c, d ) = (1)$, as claimed.

For point 3, we multiply out the norm, to get $\Norm( \alpha ) = z^p$. So 
$\alpha | z^p$. $z^p= \prod_{c \in P^*} \sigma_c(\alpha)$, and thus 
$z^p/\alpha = \prod_{c \in P^*, \sigma_c \neq Id} \sigma_c(\alpha)$. 
Consequently $\left( \alpha, z^p/\alpha \right) = (1)$, by point B. 
For the characteristic ideal, this implies: 
\[ \eu{A}^p = \left( \alpha^p, \alpha^{p-1} z, \ldots, \alpha z^{p-1}, 
\alpha \cdot (z^p/\alpha) \right) = ( \alpha ) \cdot J , \]
where the ideal $J = ( \alpha^{p-1}, \ldots, z^{p-1}, z^p/\alpha ) = ( \alpha, z^p/\alpha ) = (1)$, 
hence $\eu{A}^p = (\alpha)$: the characteristic
ideal is either principal or it has order $p$\footnote{The order of an ideal is naturally defined as the order of its class in the class group.}. 
The relation $\Norm( \eu{A} ) = ( z )$ follows from $\Norm( \alpha ) = z^p$. 

For point 4, \rf{al'} follows from $\frac{\alpha}{\bar{\alpha}} = \left( \frac{1-\bar{\zeta}}{1-\zeta} \right)^e . \frac{1+\zeta \frac{x}{y}}{1+\bar{\zeta} \frac{x}{y}}$ and $\frac{1-\bar{\zeta}}{1-\zeta} = - \bar{\zeta}$. For
\rf{twist} we note that for $c \in P^*$, 
\[ \zeta^c = (1 - \lambda )^c = 1 - c\lambda + O( \lambda^2 ), \]
and $y + \zeta x \equiv (x+y) - \lambda x \equiv (x+y) \cdot \zeta^{x/(x+y)} \bmod \lambda^2$. 

\begin{comment}
    
We define: 
\begin{eqnarray}
\label{etadef}
\eta = \upsilon \cdot \zeta^{-(1-e) x/(x+y)}
\end{eqnarray}
\end{comment}
\end{proof}
\subsection{The Stickelberger ideal and its action}
The Stickelberger element $\vartheta = \frac{1}{p} \sum_{c=1}^{p-1} c \sigma_c^{-1} \in \frac{1}{p} \Z[ G ]$ generates
the Stickelberger ideal  in the group ring of $G$ over the rational integers, by intersecting its principal ideal with $\Z[G ]$, 
according to 
\begin{eqnarray}
\label{stickid}
 I = \vartheta \Z[ G ] \cap \Z[ G ]. 
 \end{eqnarray}
 
 Comparing to the definition of the orthogonal idempotents in \rf{idemps}, we note that
 \begin{eqnarray}
 \label{st-el}
 \vartheta &=& \frac{p-1}{p} ( \veps_1 - A p ), \quad A \in \Z_p[ G ].
 \end{eqnarray}
 From the definition, $(1+\jmath) \vartheta = \Norm_{\K/\Q}$, and thus
 \begin{eqnarray}
     \label{pmstick}
     I^+ := ( 1+ \jmath ) I = \Norm \cdot I, \quad I^- := (1-\jmath) I \cong I/( \Norm). 
 \end{eqnarray}
 \subsubsection{Generators and relations}
The ideal $I$ has the property of annihilating the class group of $\K$ ( \cite{Wa}, \S 15.1). 
That is, for each ideal $\eu{C} \subset \Z[ \zeta ]$ and each $\theta \in I$, 
the ideal $ \eu{C}^{\theta} \subset  \Z[ \zeta ]$ is principal. There exists a base for 
$I^-$, made of $(p-1)/2$ elements, called {\em Fueter elements}, e.g. \cite{Mi2}, 
which are
\begin{eqnarray}
\label{fueter} 
& & \psi_n =  \vartheta ( 1+\sigma_n -\sigma_{n+1}) = \sum_{c \in S_n } n_c \sigma^{-1}_c  \in \Z_{\geq 0}[ G ], \quad \hbox{for } n \in \left\{ 1, 2, \ldots, \frac{p-1}{2} \right\}\\ 
\label{Fueter1or0}
& & \hbox{with} \quad \nonumber n_c = \left( \left[ \frac{(n+1)c}{p} \right] - \left[  \frac{nc}{p} \right] \right), \quad \hbox{and} \quad  n_c + n_{p-c} = 1, \\
%% & & \label{Fueter1or0} \hbox{and} \quad  n_c + n_{p-c} = 1,
\end{eqnarray}
where the support $S_n \subset \{1, 2 \ldots, p-1 \}$, satisfies\footnote{The set 
$p-S_n$ designates naturally $\{ p - r \ : \ r \in S_n \} $.}  
$S_n \cup ( p - S_n ) = P^*$ and
is deduced from the definition of $\psi_n$. We note that \rf{Fueter1or0} implies 
that, for $\psi_n = \sum_{c \in S_n } n_c \sigma^{-1}_c$ an element of the Fueter 
base, either $n_c = 0$ or $n_c=1$, and also that $(1+\jmath) \cdot \psi_n = 
\Norm_{\Q(\zeta)/\Q}$, in conformity with \rf{pmstick}. We shall denote conjugates 
$\psi = \sigma \psi_n$ also by {\em Fueter elements}, so in our notation, a Fueter 
element is
an element $\psi = \sum_{c \in P^*} n_c \sigma_c^{-1}$ with $n_c \geq 0$ and $n_c + n_{p-c} = 1$; in particular, $\psi + \jmath \psi = \Norm$.

We can compute, for instance, $\psi_1$:
\begin{eqnarray}
\label{psi2}
\psi_1 = \sum_{c > p/2} \sigma_{c}^{-1} \quad \hbox{and} \quad \jmath \psi_1= \sum_{c < p/2} \sigma_{c}^{-1},
\end{eqnarray} 
so $S_1 = \{ \frac{p-1}{2}, \ldots, p-1 \}$ and $w(\psi_1) = w(\jmath \psi_1) = \frac{p-1}{2}$. 

\begin{comment}
Indeed,
\begin{eqnarray}
\eu{C}^{(1+\jmath).\psi_n} =\eu{C}^{\sum_{S_n} n_c \sigma^{-1}_c + \sum_{p-S_n} n_c \sigma^{-1}_c} =  \eu{C}^{\sum_{c \in P^*} \sigma^{-1}_c} = \prod_{c \in P*} \eu{C}^{\sigma^{-1}_c} = \prod_{c \in P*} \eu{C}^{\sigma_c} = \Norm(\eu{C})
\end{eqnarray}
\end{comment}

We can write any $\theta \in I$ as
\begin{eqnarray}
\label{fueterbase} 
 \theta = \sum_{n=1}^{(p-1)/2} \nu_n \psi_n = \sum_{c=1}^{p-1} n_c \sigma_c^{-1}; \quad \nu_n, n_c \in \Z. 
 \end{eqnarray}

Therefore, for each ideal $\eu{C} \subset \Z[ \zeta ]$ and each $\theta \in I$, the ideal $ \eu{C}^{\theta} \subset  \Z[ \zeta ]$ is generated by some $\gamma \in \Z[\zeta ]$, which satisfies $\gamma \cdot \overline{\gamma} = \Norm( \eu{C} )^{\varsigma (\theta)}$, for an integer $\varsigma (\theta) \in \Z$, which we call the {\em relative weight} of $\theta$.
The {\em absolute weight} (or simply, {\em weight},) of  $\theta = \sum_{c \in P^*} n_c \sigma^{-1}_c \in \Z[ G ]$ is $w( \theta ) = \sum_c | n_c |$. We say that $\theta = \sum_{c \in P^*} n_c \sigma_c^{-1}$ is {\em positive}, writing $\theta \in I_+$, if $n_c \in \Z_{\geq 0}$ for all $c \in P^*$. Note that the relative weight of each Fueter element is $1$.

\subsubsection{The Fermat quotient ideal, $J_k \subset I$ and $I/pI$}
We define the Fermat quotient map $\phi : \Z[ G ] \ra \F_p$ by
$  \theta \mapsto \phi( \theta ) \quad : \quad \zeta^{\theta} = \zeta^{\phi( \theta)}$. Explicitly,
\begin{eqnarray}
\label{phi}
\phi \left(  \sum_{c \in P^*} n_c \sigma_c^{-1} \right) = \sum_{c \in P^*} n_c/c \in \F_p. 
\end{eqnarray}
We identify the value $\phi( \theta ) \in \F_p$ with its natural lift to $\N$,
under the least positive remainder representation of $\F_p$. 
\begin{definition}
\label{FID}
    The {\em Fermat ideal} is $I_0 = I \cap \Ker( \phi )$:
this is the module of all Stickelberger elements $\theta$ such that $\zeta^{\theta} = 1$.
The module $J_k \subset I_0$ is defined by
\begin{eqnarray}
    \label{jk}
    J_k = \{ \theta \in I_{0,+} \ : \ \varsigma( \theta ) = 2k, \ k \geq 1 \}.
\end{eqnarray}
It is thus the submodule of $I_0$ consisting of elements that are sums $\theta = \theta_1 + \theta_2$
with $\varsigma( \theta_i ) = k; i = 1, 2$; the $\theta_i$ need not be elements of $I_0$.
\end{definition}
The following fact shows that we always can choose elements $\theta \in I_{0,+}$
of small relative weight:
\begin{fact}
 \label{i0}
 For $p \geq 5$ there always exists an element $\theta \in I_{0,+}$ with 
 $\varsigma (\theta)\geq 2$. In particular, $J_k \neq \emptyset$ for all $p > 3$
 and $k \leq \frac{p-1}{2}$
\end{fact}
\begin{proof}
 Let $\phi( \psi_1 ) = a$ and $\phi( \psi_2 ) = b$. 
 If $a \cdot b \equiv 0 \bmod p$, then there exists $j \in \{ 1, 2 \}$ 
 such that $\phi( \psi_j ) = 0$, so $\theta = 2 \psi_j$
 satisfies the claim. Otherwise, let $c \in \N$ be such that
 $ a + b \cdot c \equiv 0 \bmod p$. Since $\zeta^{\phi( \theta ) } = \zeta^{\theta}$, 
 it follows that $\phi( \sigma_c \theta ) \equiv c \phi( \theta ) \bmod p$, 
 and thus $\phi( \psi_1 + \sigma_c \psi_2 ) = 0$. 
 Therefore, $\theta = \psi_1 + \sigma_c \psi_2$ satisfies the claim. 
 Since $I^-$ contains, for $p \geq 5$, at least two $\Z$-base elements, the claim follows.
 The same procedure can be applied a fortiori for 
 larger relative weights, thus confirming the claim on $J_k$. 
\end{proof}

The Stickelberger ideal annihilates in particular the characteristic ideal 
$\eu{A}$. By Lemma \ref{aux1}, the ideal $( p, \Norm ) \Z[ G ] \subset \Z[ G ]$ is 
a natural annihilator of $\eu{A}$; we are interested in {\em non-trivial annihilators} in $I$, and this will be such $\theta \in I$ which have 
non trivial images in 
\begin{eqnarray}
    \label{Ip}
    I_p := I^-/p I^- = I/( p, \Norm )
\end{eqnarray}, 
according to \rf{pmstick}. The spectral decomposition \rf{idemps}
induces 
\[ I_p = \sum_{c \in P^*; c \equiv 1 \bmod 2} \veps_c I_p, \]
in which $\veps_c I_p = 0$ if and only if the Bernoulli number $B_{p-c} \equiv 0 
\bmod p$, or, equivalently, $\veps_c \id{C}( \K )[ p ] \neq \{ 1 \}$  -- 
see \cite{Wa}, \S 6.3 for details. 
\begin{remark}
    \label{e1}
    We note the following fact for future reference: suppose that 
    $0 \neq t \in J_1$
    and $a = \sigma - \varpi( \sigma ) \in \F_p[ G ]$. Obviously, $t$ has non
    trivial image in $I_p$, since its coefficients are $0 \leq n_c \leq 2$ and
    they cannot all vanish modulo $p$. Note also that 
    $\veps_1 I_p = \F_p \cdot ( p \vartheta)$. There is thus some 
    $c \in P^* \setminus \{ 1 \}; c \equiv 1 \bmod 2$ such that 
    $\veps_c t$ has non trivial image in $I_p$. Since $\sigma \veps_c = \varpi( \sigma )^c \veps_c$, it follows that 
    \begin{eqnarray}
        \label{nonv}
        \veps_c ( a t ) = (\varpi( \sigma )^c - \varpi( \sigma )) \veps_c \not \equiv 0 \bmod (p, \Norm).
    \end{eqnarray}
\end{remark}

\subsubsection{The action of the Stickelberger ideal on characteristic ideals}
\label{sti}
Since $\theta \in I$ annihilates the class group, there is some principal ideal $b(\theta)$ such that
$b(\theta) = \eu{A}^{\theta}$ and it satisfies $b(\theta) \cdot \bar{b}(\theta) = \Norm( \eu{A} ) = ( z )$.
It is known from the theory of Gauss and Jacobi sums  -- see e.g.
\cite{Wa} -- that principal ideals arising from the action 
of the Stickelberger ideal are generated by 
{\em Jacobi numbers}, which are defined as products of Jacobi sums.

Iwasawa proved in \cite{Iw} that Jacobi numbers $\rg{J}$ verify\footnote{In the classical definition $\tau(\chi) = \sum_{x \in \ZMs{q}} \chi( x ) \xi^x$,
Iwasawa actually proves that $\tau(\chi) \equiv -1 \bmod \eu{P}$, with $\eu{P} \in \Z[ \zeta_p, \xi ]$ 
an ideal above $p$. This led Lang to modify the definition by changing the sign, 
and we use his definition here.}

\begin{eqnarray}
\label{Iw}
 \rg{J} \equiv 1 \bmod ( 1 - \zeta )^2,
\end{eqnarray} 
Since the product of Jacobi numbers by their complex conjugates are rational integers, the above
condition implies that there is a unique Jacobi number $\beta$ that generates the ideal $b$, and all
other generators of this ideal, which are rational upon multiplication by their complex conjugates, 
differ from the Jacobi number $\beta$ by a root of unity -- a consequence
of the Kronecker unit theorem. 
\begin{comment}
Let $\beta(\theta) \in b(\theta)$ be the unique Jacobi number generating the ideal $b(\theta)$. Since 
$(\alpha^{\theta}) = \eu{A}^{p \theta} = b(\theta)^p$, we obtain by using Lemma \ref{aux1}, point 4, 
in combination with \rf{Iw}, that
\begin{eqnarray}
\label{upto}
\alpha^{(1-\jmath) \theta} = \eta' \cdot {\beta(\theta)}^p, \quad
\eta' = \begin{cases}
         \zeta^{-2x \theta/(x+y)} & \hbox{ for $e = 0$, and}\\
         (-\zeta^{-\theta}) & \hbox{otherwise.}
        \end{cases}
\end{eqnarray}
For $e = 1$, we have to choose a Fueter element in \rf{fueter} which satisfies $\zeta^{\vartheta} \neq 1$. Because
$\phi( \psi_n ) + \phi( \psi_{p-n} ) \equiv -1 \bmod p$, 
for each $n$, either $\psi_n$ or $\jmath \psi_n$ will satisfy our condition. 
\end{comment}

\subsubsection{The $\beta$-map}
\label{betamap}
For $\theta \in I$, let the auxiliary number $ \beta( \theta )$ 
be the Jacobi number generating $\eu{A}^{\theta}$.
\begin{lemma}
    \label{betahom}
    The map 
    \begin{eqnarray}
    \label{bmap}
        \beta: I \ra \K^{\times}, \ \ \theta \mapsto \beta(\theta) \nonumber
\end{eqnarray} defines an injective homomorphism of $G$-groups, via 
    \begin{eqnarray}
    \label{bhom}
         \beta( \theta_1 + \theta_2 ) = \beta( \theta_1 ) \cdot \beta(\theta_2 ); \quad \beta( - \theta ) = 1/\beta( \theta ).
    \end{eqnarray}
\end{lemma}
\begin{proof}
    The homomorphism relations in \rf{bhom} can easily be verified from the definition. 
    For $\sigma \in G$, we have, by definition 
    \[ \beta( \sigma \theta ) = \beta( \theta )^{\sigma} \in \K^{\times},\]
    so the homomorphism is one of $G$-groups. For injectivity, 
    suppose that $\theta = \sum_c n_c \sigma^{-1}_c $ and $\beta( \theta ) = 1$. Then 
    \begin{eqnarray*}
        \beta( \theta )^p  =  {\alpha'}^{\theta} = \prod_{c=1}^{p-1} \sigma_c^{-1}( \alpha' )^{n_c} = 1.
    \end{eqnarray*}
    By Lemma \ref{aux1}, the ideals $\left( \sigma_c^{-1}( \alpha' ) \right)$ are pairwise 
    coprime, so the exponents $n_c$ must all cancel.
\end{proof}
\begin{remark} 
     \label{furt}
     We can now derive Furtw\"angler's result. Suppose that $(x, y, z)$ is a solution to \rf{flt} with $p \nmid x y z$. Then $x/(x+y) \equiv d [p]$ for some $d \in P^*$
     and, if $t \in I$ is such that $\zeta^t \neq 1$, then \rf{twist} implies that 
     \[ \zeta^{-d t} \equiv \beta( t )^p \bmod x \Z[ \zeta ].\]
     Thus, the equation $\zeta^{d'} := \zeta^{-d t} \equiv T^p \bmod \eu{R}$ has the non trivial 
     solution $T = \beta( t ) \in \Z[ \zeta ]/\eu{R}$
     for each prime ideal $\eu{R} \subset \Z[ \zeta ]$ with $x \equiv 0 \bmod \eu{R}$. Therefore, 
     these ideals split in $\Q[ \zeta_{p^2} ]/\K$. By interchanging the roles of $x, y, z$, we obtain Furtw\"angler's result \cite{Fur}.
 \end{remark}

\subsubsection{Some properties of the module $J_k$}
\label{sjk}
We start by noticing that for  all $t \in J_k$, we have ${\alpha'}^t = \alpha^t$ and thus
\begin{eqnarray}
\label{pt}
\beta( t )^p = \alpha^t. 
\end{eqnarray}
This is a consequence of $\zeta^t = 1$. 
We note the following actions of $t \in J_k$ on $\lambda$:
\begin{eqnarray}
\label{pm}
    ( 1 - \zeta )^{t( 1+ \jmath)} & = & p^{2k}; \quad ( 1 - \zeta )^{t( 1 - \jmath)} = (-\zeta)^t = (-1)^{k(p-1)} \cdot \zeta^t = 1, \quad \hbox{hence} \
    \lambda^t = \pm p^k
\end{eqnarray}

\begin{definition}
    \label{alredef}
    We shall, from now on, only work with $t \in J_k$. In view of \rf{pt} and \rf{pm}, and the fact that we will only be using relation \rf{al'} from Lemma \ref{aux1}, we redefine $\alpha = y + \zeta x $ also for the case $e = 1$. In this case, this redefined $\alpha$ verifies the previous relation and allows for a unified treatment of both cases of \rf{sft}, both when $( K, p ) = 1$ and when $p | K$.
\end{definition}
In the case $e = 1$, the first line in \rf{pm} implies that
$\beta(t)^p = s(t) \alpha^{t}/p^{k}$ and thus
\begin{eqnarray}
\label{betasm1}
\beta(t)^{\sigma-1} = \alpha^{t( \sigma-1 )} , \quad 
 \forall \sigma \in G \setminus \{1\}.
\end{eqnarray}
showing that the denominator vanishes in the case $e = 1$, when acting with $\sigma - 1$ on $\beta( t )$,
for arbitrary $\sigma \in G \setminus \{ 1 \}$ and $t \in J_k$.
The following result induces the key factorizations of $K$, mentioned in the introduction.
\begin{lemma}
    \label{psipol}
    Let $\Psi_t( x, y ) = \alpha^t$ and remember that $K = x y (x+y )$. Then 
    $K | \left( \Psi_t (x, y ) - \sigma(\Psi_t (x, y ))\right)$ for all $\sigma \in G  \setminus \{1\}$.
\end{lemma}
\begin{proof}
    We have the following congruences:
    \begin{eqnarray*}
        \Psi_t( x, y ) & \equiv & \begin{cases}
            y^{(p-1) k} \bmod x, \\
            x^{(p-1) k} \bmod y, \\
            \lambda^t y^{(p-1)k} = s(t) p^k y^{(p-1) k} \bmod [x+y], s(t) \in \{-1,1\}.
        \end{cases}
    \end{eqnarray*}
    The right side of the congruence being rational in all cases, the claim follows.
\end{proof}

\section{The Strong Fermat Equation}
We assume that \rf{sft} has a non trivial solution. With our redefined $\alpha$, $\Norm( \alpha ) = p^e \cdot z^p$.
Let $t \in J_k$. Let $\beta( t )$ be the Jacobi sum generating $\eu{A}^t$.
Since $\beta(t) \equiv \sigma \beta(t) \bmod \lambda$, it 
follows that $\beta(t)^p - \bar{\beta}(t)^p$ is an algebraic 
integer, divisible by $p$. In view of Lemma \ref{psipol}, it is also divisible by $K$; in
the sequel we shall focus on the non $p$-part of $K$ and define the following:
\begin{definition}
\label{datadef}
Let $t \in J_k$, $w \in F = \{ x, y, x+y\}$ and $' : F \ra \Z$ the map $v \mapsto v/p^{v_p( v )}$. For each pair $(t, w)$ and $\sigma \in G \setminus \{ 1 \}$, let:
\label{data1}
 \begin{eqnarray*}
     \Delta( \sigma ) & = & \Delta_{t, w}( \sigma )  =  \Psi_t - \sigma \Psi_t, \\
     \delta_c(\sigma ) & = & \beta(t) - \zeta^c \sigma(\beta(t)), \quad c \in P, \\
     \eu{D}_c(\sigma ) & = & ( w', \delta_c( t, \sigma ) ), \quad c \in P, \\
     \rg{D} & = & \bigcap_{\sigma \in G} \eu{D}_0( \sigma ).
 \end{eqnarray*}    
\end{definition}
The items in the above definition depend on a set of variables: $t, w, \sigma, c$. In the 
sequel, we shall always keep the reference on $c$ explicit, while the further variables
will only be mentioned when their choice changes or is not obvious in the context. By definition, 
$(w', p ) = 1$ and a fortiori, $( \eu{D}_c( \sigma ), p ) = 1$ for all $\sigma \in G$.

As an immediate consequence of the definitions and of Lemma \ref{aux1}, we have
\begin{lemma}
    \label{cop}
    Let $t \in J_k; w \in F = \{ x, y, x+y\}$ and $\sigma \in G \setminus \{ 1 \}$ be fixed.
    Then
    \begin{itemize}
        \item[ 1. ] The ideals $\eu{D}_a, \eu{D}_b$ are pairwise coprime  
        for distinct $a, b \in P$, and real for $a \in P$ and $\sigma = \jmath$.
    \item[ 2. ] The product  $\prod_{a \in P} \eu{D}_a = ( w' )$. 
    \item[ 3. ] The ideal $\rg{D}$ is $G$-invariant, so $\rg{D} = (D)$ for some $D \in \N$.
    Moreover, $\beta( t ) \equiv \sigma( \beta ( t )) \equiv z^k \bmod \rg{D}$
    for all $\sigma \in G$.
    \item[ 4. ] Let $\eu{R} \subset \Z[ \zeta ]$ be a prime that divides $w'$. If $\eu{R} \nmid \rg{D}$, then there is a $\sigma \in G$ and $a \in P^*$ such that $\eu{R} | \eu{D}_a( \sigma )$.
    \end{itemize} 
\end{lemma}
\begin{proof}
    Defining the ideal $D( a, b ) = ( \delta_a, \delta_b )$, we see that for $a, b > 0$:
    \begin{eqnarray*}
        (\bar{\sigma \zeta}^a - 1) \delta_a - (\bar{\sigma \zeta}^b - 1) \delta_b & = & (\zeta^{-a \sigma} - \zeta^{-b \sigma}) \beta(t) \in D( a, b ), \\
        (1-\zeta^a) \delta_a - ( 1-\zeta^b )\delta_b & = & ( \zeta^b - \zeta^a ) \sigma(\beta(t)) \in D( a, b ), \quad \hbox{hence} \\
        D( a, b ) & \mid & ( \beta(t), \sigma( \beta(t) )) \cdot (\lambda).
    \end{eqnarray*}
    Now, $\beta(t) | z^{\varsigma( t )}$ and, since $( z, w ) = 1$, it follows that 
    $( \eu{D}_a, \eu{D}_b ) | ( \lambda )$.
    Since $( \eu{D}_a, p ) = 1$, we conclude that $( \eu{D}_a, \eu{D}_b ) = 1$.
    A similar argument holds when one of $(a, b )$ is $0$ -- one may assume $a = 0$ 
    and build appropriate differences to show that $( \eu{D}_0, \eu{D}_b ) = 1$ 
    in this case too; we leave the details to the reader.

    Note that
    \begin{eqnarray*}
        \bar{\delta}_a = (  \bar{\zeta}^a ) \cdot \delta_a \cdot \frac{z^{\varsigma( t )}}{\beta(t) \sigma( \beta(t) )} = 
        \delta_a \cdot \frac{\bar{\zeta}^a \bar{\beta}(t)}{\sigma( \beta(t) )}.
    \end{eqnarray*} 
    The complex conjugates of $\eu{D}_a( \sigma )$ are thus equal to $\eu{D}_a$ as ideals, 
    if and only if $\sigma = \jmath$, thus completing the proof of 1.
    
    We have shown that $\Psi_t(x, y ) \equiv \sigma( \Psi_t( x, y )) \bmod w$; 
    since 
    \[ \prod_{a \in P} \delta_a = \frac{\Psi_t( x.y )-\sigma( \Psi_t( x, y ))}{\prod_{a \in P} \nu_a}, \]
    and the denominators $\nu_c$ are powers of $(\lambda)$, while $( \eu{D}_a, p ) = 1$, 
    it follows that $\prod_{a \in P} \eu{D}_a = (w')$, which confirms 2.

    From the definition, if a prime $\eu{R} | \rg{D}$, then 
    $\beta(t) \equiv \sigma( \beta(t) ) \bmod \eu{R}$ for 
    all $\sigma \in G$. Thus $\beta(t) \equiv C = \frac{\Tr( \beta(t) ) }{p-1} \bmod \eu{R}$, 
    with $C \in \Q$.
    Since $\beta(t) \equiv \bar{\beta}(t) \bmod \eu{R}$, it also follows by multiplying the 
    two congruences, that $z^{2k} \equiv C^2 \bmod \eu{R}$, hence 
    $C \equiv s z^k \bmod \eu{R}$, for $s \in \{-1,1\}$. 
    We claim that the sign of $C$ is in fact positive. For some $w \in \{ x, y, x + y, x-y\}$, 
    and $\eu{R} \ | \ \rg{D} \ | \ w'$, we have, by raising to the odd power $p$,
    \begin{eqnarray*}
        \beta(t)( \sigma, w ) & \equiv & s z^k \bmod \eu{R} \quad \Rightarrow \quad \alpha^{\sigma t} \equiv s z^{pk} \bmod \eu{R}.
    \end{eqnarray*}
    Using Lemma \ref{psipol}, we note that $k \cdot \Norm \in J_k$ and the proof of the Lemma also applies to 
    $z^{pk} = \Norm( \alpha^{k \sigma })$. Consequently, $ \alpha^{\sigma t} \equiv z^{p k} \bmod K$, and $s = 1$,
    as claimed.

    By acting with $G$ on the 
    congruence we derived, we conclude that:
    \begin{eqnarray}
        \label{valD}
        \sigma( \beta(t) ) \equiv z^k \bmod \tau( \eu{R} ), \quad \forall \sigma, \tau  \in G.
    \end{eqnarray}
    The congruences thus hold modulo $r$, the rational prime below $\eu{R}$. 
    This holds for all $\eu{R} | \rg{D}$, so $\rg{D}$ is $G$-invariant. This confirms 3.
    
    Moreover, by definition of $\rg{D}$, any prime $\eu{R} | ( w' )$, coprime to $\rg{D}$
    will divide  
    $\eu{D}_a( \sigma )$ for some $a \in P^*$ and $\sigma \in G \setminus \{ 1 \}$. 
    This confirms point 4.
\end{proof}
The above Lemma indicates that the $\eu{D}_a$ induced uncommonly high factorizations of
the numbers $w \in F =\{ x,y, x+y\}$. This observation is not particularly new,
and in itself it cannot bring final insights about the equations under investigation.
However, the definition of the ideal $\rg{D}$ leads farther, and we deduce from the above:
\begin{lemma}
    \label{gcd0}
    Under the notations of Definition \ref{datadef}, we have $D = 1$.
\end{lemma}
\begin{proof}
    Let's assume $D \neq 1$. By Fact \ref{i0}, we can choose $t \in J_k$ for any $p > 3$ and $k \leq \frac{p-1}{2}$.
    We may write $t = t_1 + t_2$, with $t_i \in J_{k}; i = 1, 2$. Then
    $t = t_1 + t_2 = t_1 + k \Norm - \bar{t}_2$.
    
    For any $w'$ in Definition \ref{datadef}, and any $\sigma \in G \setminus \{ 1 \}$, we know 
    from \rf{valD} that $\sigma( \beta( t )) \equiv z^k \bmod D$. 
    There is thus a $\chi \in \Z[ \zeta ]$ such that
    \[ \beta( t ) -  z^k = D \cdot \chi. \]
    Write now $\beta( t ) = \beta( t_1 ) \cdot \beta( t_2 )$ and 
    $z^k = \beta( t_2 )^{1+\jmath}$. We then have
    \[ \beta( t_1 ) \cdot \left( \beta( t_2) - \bar{\beta}( t_1 )\right) = D \cdot \chi. \]
    Since $D | w'$ and $\beta(t) | z^{2k}$, it follows that $( \beta( t_1 ), D ) = 1$, 
    so we conclude that $\beta( t_1 ) | \chi$. We repeat the same argument, by 
    interchanging the roles of $t_1$ and $t_2$, finding that $\beta( t_2 ) | \chi$ 
    too. There is thus a $\chi' \in \Z[ \zeta ]$, such that 
    $\chi = \beta( t ) \chi'$. The initial identity becomes
    \[ \beta( t ) \cdot ( 1 - D \chi' ) = z^k .\]
    Multiplying by complex conjugates in the same identity, we
    conclude that $( 1 - D \chi' )^{1+\jmath} = 1$, and thus 
    $\mu = 1 - D \chi' \in \lan \pm \zeta \ran$, by the Kronecker Unit Theorem. 
    So $\beta( t ) = \mu^{-1} z^k$ and $\mu^{-2p} = 1$.
    But then $( \mu^{-1} z^k )^{2p} \in \Z$, while $\beta( t )^{2p} = \alpha^{2t} \not \in \Z$,
    as follows from Lemma \ref{aux1}. 
    The assumption $D \neq 1$ is thus untenable; of course, if $D = 1$, 
    the original congruences modulo $D$ are void, so there is no contradiction.
\end{proof}

We came a step further, since the Lemma implies that for any $t \in J_k$, 
every rational factor of $w'$
will have a non trivial factor in some ideal $\eu{D}_a( \tau )$, 
for some $\tau \in G \setminus \{ 1 \}$. This fact, together with the result of Lemma \ref{gcd0}
allows us to show that these factors split completely in $\K^+$:
\begin{lemma}
    \label{gcd1}
    Let $w' \in F' = \{ x', y', (x+y)' \}$, as per Definition \ref{datadef}.
    Then every prime $r | w'$ is totally split in $\K^+$.
\end{lemma}
\begin{proof}
Let $t = t_1 + t_2 \in J_1$ and consider a prime $r | w'$. Assuming that $r$ is not
    totally split in $\K^+$, we let $\eu{R} | r$ be a prime above $r$ in $\K$. 
    Let $D_{\eu{R}} \subset G$ be the decomposition group of $\eu{R}$. 
    For any $\tau \in G$, there is an $a( \tau ) \in P$ such that 
    $\eu{R} | \eu{D}_{a( \tau )}( \tau )$.
    Thus
    \begin{equation}
    \label{congrus}
        \beta( t )^{(1-\tau) }  \equiv \zeta^{a( \tau )} \bmod \eu{R},  \ \forall \tau \in G.
    \end{equation}
    Let $\nu \in G$ be such that $\jmath \nu \in D_{\eu{R}}$ fixes $\eu{R}$. For any
    such $\nu$, raising \rf{congrus} to $\varpi( \nu ) + \nu$ is an operation that
    fixes $\eu{R}$, yielding
    \begin{equation}
    \label{congruclean}
        \beta( t )^{(\varpi( \nu ) + \nu)(1-\tau) }  \equiv \zeta^{a( \tau )(\varpi( \nu ) + \nu)} \equiv 1 \bmod \eu{R},  \ \forall \tau \in G.
    \end{equation}
    
    By Remark \rf{e1} we may choose $\nu$ such that $\theta = (\varpi( \nu ) + \nu) t \in J_{k} \setminus (p , \Norm ) I_+$. By \rf{congruclean}, for such $\theta$ we see that 
    $\eu{R} | \rg{D}( \theta )$. Since the right hand side ideal $\rg{D}( \theta ) = 1$
    by Lemma \ref{gcd0}, it follows that $\eu{R}$ must be trivial.
    Plainly, there are no prime divisors of $w'$ that are not totally decomposed in $\K^+$, which completes the proof.
\end{proof}

We thus conclude that
\begin{proposition}
    \label{pmain}
    Let 
    \[ K' = \frac{K}{p^{v_p \left( K \right)}}. \]
    The primes $r | K'$ are all totally split in $\K^+$, so $r \equiv \pm 1 \bmod p$.
    Moreover, if $r | K$ then, either $r | K'$ or $r = p$.
\end{proposition}
\begin{proof}
    Let $F = \{ x, y, x+y\}$. Note that $K$ and $K'$ only differ by a power of $p$. 
    By Lemma \ref{gcd1}, the primes $r | K'$ are totally split in $\K^+$, 
    thus verifying $r \equiv \pm 1 \bmod p$. The second claim is an obvious consequence of
    the definition of $K'$.
\end{proof}
As a consequence:
\begin{corollary}
    \label{thm}
    Equation \rf{sft} has no solutions for $p > 3$ and Theorem \ref{main} is true.
\end{corollary}
\begin{proof}
    At least one $v \in F$ must be even: indeed, if $x$ and $y$ have the same parity, 
    then $x + y \equiv 0 \bmod 2$. Otherwise, one of $x$ or $y$ must be even, 
    thus confirming the claim. Since $p$ is odd, $v$ and $v'$ have 
    the same parity for all $v \in F$. Since $p$ is odd, it remains that $2 | K'$. By Proposition
    \ref{pmain}, all the primes $r | K$
    are of the form $r = m p \pm 1$ or $r = p$, so this should hold also for 
    $r = 2$. But $p > 3$, so we reach a contradiction, which confirms Theorem \ref{main}.
\end{proof}

\begin{remark}
\label{Wiles}
We obtained along with the general result, a new, cyclotomic proof of the popular Fermat's Last Theorem. Of course, 
after the epochal proof of Wiles \cite{W} and Wiles and Taylor \cite{WT}, this is just an alternative approach -- 
having the advantage to allow proof of the non-existence of solutions to a  wider range of results, for instance, 
equations of type \rf{flt} in abelian fields, or twisted by a constant on the right hand side.

%%It is important to mention at this place that another alternative proof of FLT, which is based, like \cite{W}, on a 
%%major theoretical breakthrough was recently achieved: after the publishing of the Inter Universal Teichm\"uller 
%%Theory \cite{Mo}, Mochizuki achieved together with four other authors, an effective version of the 
%%$abc$-inequality, and the result is accepted for publication: the proof comes with yet an other proof of FLT 
%%\cite{MFHMP}. The present proof arouse while completing the details for some improved lower bounds needed in 
%%\cite{MFHMP}.
\end{remark}

\vspace*{2.0cm}
\textbf{Acknowledgements}
We thank the referee who contributed with their questions and remarks to the improvement of the final version of the present paper. 
This work is built on shoulders of giants, mainly on the theories and results due to Kummer, Hensel and Chevalley. 
One cannot however under-appreciate the contribution of numerous later authors, whose work with these methods made a combination of ideas which might have seemed remote in earlier times, become possible, in a natural way.

We are deeply indebted to the people who contributed to this recent surprising development. We are grateful to  I. Fesenko and the students H. Chen, D. Lehmann, C. Liu and J. Reichardt for their active interest and questions which helped improve the presentation. 

It is not possible to list, or even
recall, the sporadic yet important discussions and challenges which made that, despite the decision to stay away from FLT, the interest in the question was warmed up periodically, so it never completely disappeared, over more than two decades. The  support of the colleagues at the Mathematical Institute of the University of G\"ottingen will be mentioned, collectively -- to all, our thanks.

Most of all, the second author is grateful to the family in small -- Theres and Seraina -- and at large, unlistable, who stood beside during short and 
long times.  The first author is grateful to his close family -- Corinne and Emma -- who have supported him through ups and downs, and to Y. Bilu and P. Mih\u{a}ilescu, who have allowed him to come back onto a path he should have never left.

\end{document}